\documentclass[11pt,a4paper]{article}
\usepackage{color,amsmath,amssymb, lscape,graphicx,ifthen,pdfpages}
\usepackage[colorlinks=true,
linkcolor=webgreen,
filecolor=webbrown,
citecolor=webgreen]{hyperref}
\definecolor{webgreen}{rgb}{0,.5,0}\definecolor{webbrown}{rgb}{.6,0,0}

\newcommand{\seqnum}[1]{\href{http://oeis.org/#1}{\underline{#1}}}

\def\ogf{o.g.f.\,}

\def\eg{{\it e.g.},\,}

\def\ie{{\it i.e.},\,}

\def\via{{\it via}\, }
\def\sspp{\,+\,}
\def\sspm{\,-\,}
\def\sspeq{\,=\,}
\def\sspdef{\, :=\,}
\def\sspneq{\,\neq\,}
\def\sspkl{\,<\,}

\def\sspgeq{\,\geq\, }
\def\sspleq{\,\leq\,} 
\def\pn{\par\noindent}
\def\pb{\par\bigskip}

\def\pbn{\par\bigskip\noindent}
\def\psn{\par\smallskip\noindent}

\def\Beq{\begin{equation}}
\def\Eeq{\end{equation}}
\def\Beqarray{\begin{eqnarray}}
\def\Eeqarray{\end{eqnarray}}
\def\Eq#1{Eq.\,#1}

\def\sspin{\,\in\,}

\def\sspto{\,\to\,}
\def\rhs{right-hand side\,}
\def\lhs{left-hand side\,}

\def\floor#1{\left\lfloor{#1}\right\rfloor}

\def\binomial#1#2{{#1} \choose {#2}}
\def\range#1#2{#1,\,\count15=#1 \advance\count15 by +1 \number\count15,\,...,\,#2\,}
\def\rangeinf#1{#1,\, \count16=#1 \advance\count16 by +1 \number\count16,\,...\,}

\def\LTmatz#1#2#3#4{\left(\begin{array}{cc}#1&#2 \\ #3&#4 \end{array}\right)} 
\def\LTmatd#1#2#3#4#5#6#7#8#9{\left( \begin{array}{ccc}#1&#2&#3 \\ #4&#5&#6 \\ #7&#8&#9 \end{array}\right)}  
\normalbaselines
\begin{document}
\bibliographystyle{unsrt}
\rightline{April 2023}
\vbox {\vspace{6mm}}
\begin{center}
{\Large {\bf Improved Formula for the Multi-Section of the Linear Three-Term Recurrence Sequence}}\\ [9mm]
{\bf Gary Detlefs} \footnote{\href{mailto:gdetlefs@aol.com}{\tt gdetlefs@aol.com}}{\bf and Wolfdieter Lang} \footnote{ 
\href{mailto:wolfdieter.lang@partner.kit.edu}{\tt wolfdieter.lang@partner.kit.edu},\quad 
\url{http://www.itp.kit.edu/~wl}
                                          } \\[3mm]
\end{center}
\begin{abstract}
The standard formula for the multi-section of the general linear three-term recurrence relation is simplified in terms of Chebyshev S-polynomials.
\end{abstract}
\section {\bf Introduction}

 \hskip 1cm The $m-$section (multi- or modular-section) of an integer sequence consists of set of $m$ sequences which carry as indices the equivalence classes modulo $m$. 
\pn
The general decomposition of the ordinary generating function (\ogf\llap) of the sequence into the $m$ \ogf\llap s \pn
of the members of the set of $m-$sections is given in terms of a special $m \times m$  {\sl Vandermonde} matrix. The inverse of this matrix gives the \ogf\llap s of these members in terms of the \ogf of the sequence. The computation which brings these $m$ fractions into one, either by hand (tedious) or by computer, does not give an insight into the structure of this final rational fraction.
\pn
For the general sequence satisfying a linear three-term recurrence relation (called {\sl Horadam}-sequence) it is shown that the result for the \ogf\llap s of the $m-$section sequences can be given in terms of {\sl Chebyshev}-$S$ (and -$R$ polynomials, the monic $T$-polynomials, which are a difference of two $S$-polynomials).
\pn
This is achieved by a proposal for the $m$-section of the {\sl Horadam} sequence, first a conjecture for the first element of this section by one of the authors ({\sl G. D.}), then generalized for all elements, and later proved by the second author.
\psn
The first section summarizes the standard treatment of the $m$-section of a sequence and the \ogf\llap s. The second section is a reminder of some elementary properties of the {\sl Horadam} sequence. In the third section the conjectures for the $m$-section of this sequence are formulated, and the last section gives the proof of these conjectures.
\pn
The proof uses a lemma a (known) alternative bisection of the {\sl Chebyshev}-$S$ polynomials (not the one obtained for improved $m\sspeq 2$ case).
\section {\bf Multi-Section of a sequence} 
\hskip .5cm This section is a reminder of the standard treatment of the $m-$section of a sequence.\pn
The ordinary generating function (\ogf\llap) $G(m,l,x)\sspeq \sum_{n=0}^{\infty} a(m\,n+l)\,x^n$ of the $l$th part of the $m$-section of a sequence $\{a(n)\}_{n>=0}$ with \ogf $G(x)\sspeq  \sum_{n=0}^{\infty} a(n)\,x^n$, for integer $m \sspgeq 2$ and $l\sspin \{\range{0}{m-1}\}$, satisfies
\Beq
G(x) \sspeq \sum_{l=0}^{m-1} G(m,l,x^m)\,x^l\,.
\Eeq
For the solution of $G(m,l,x)$ for given $G(x)$ one uses the roots of the polynomial $x^m\sspm 1$, that is $w(m,\,k)\sspeq e^{2\pi k/m}$, for $k\sspin \{\range {0}{m-1}\}$, and considers the inhomogeneous system of $m$ equations, for $k\sspin\{\range{0}{m-1}\}$, for the $m$ unknowns $\{G(m,l,x)\}_{l=0}^{m-1}$,using a {\sl Vandermonde} matrix $V_m(x)$ with elements 
\Beq
[V_m(x)]_{k,l} \sspeq (w(m,\,k)\,x)^l,
\Eeq
as
\Beq
\sum_{l=0}^{m-1} [V_m(x)]_{k,l}\,G(m,l,x^m) \sspeq G(w(m,\,k)\,x)
\Eeq
Note that $(w(m,\,k)\,x)^m \sspeq x^m$ has been used.
\pn
The inverse of a general {\sl Vandermonde} matrix is known. \eg \cite{Knuth}, and for the present case its elements become
\Beq
[V^{-1}_m(x)]_{l,j} \sspeq N(m,l,j,x)/DN(m,j,x),
\Eeq
with denominator
\Beq
DN(m,j,x) \sspeq x^{m-1}\,\prod_{0\sspleq \sspkl i\sspneq j\sspleq m-1}\, (w(m,\,i)\sspm w(m,\,j)), 
\Eeq
and numerator
\Beq
N(m,l,j,x) \sspeq (-1)^l\,x^{m-1-l}\, \sum_{n=1}^{\#Ch(m,l,j)}\,\prod_{k=1}^{m-1-l} (Ch(m,l,j)[n])[k],
\Eeq
where the list of lists (order respected, and the $k$th elements of a list $L$ is denoted by $L[k]$) 
\Beq
Ch(m,l,j)\sspeq choose(P(m,\,j),m-1-l),
\Eeq
with the list
\Beq
P(m,\,j) \sspeq [w(m,\,0),...,w(m,\,j-1), w(m,\,j+1),\, ...,\,w(m,\,m-1)]\ .
\Eeq
The length of list $Ch(m,l,j)$ is $\#Ch(m,l,j)\sspeq {\binomial{m-1}{l}}$ and the length of the lists of $Ch(m,l,j)$ is $m-1-l$ with $\#P(m,\,j)\sspeq m-1$.\pn
Thus, using new arguments $x \to x^{1/m}$, one obtains, for  $l\sspin \{\range{0}{m-1}\}$
\Beq
G(m,l,x)\sspeq \sum_{j=0}^{m-1} [V^{-1}_m(x^{1/m})]_{l,j}\,G(w(m,j)\,x^{1/m}). 
\Eeq
\psn
{\bf Example 1: $\bf m\sspeq 3$}\psn
With $w(3,\,0)\sspeq 1$, $w(3,\,1)\sspeq w \sspeq \frac{1}{2}(-1\sspp \sqrt{3}\,i)$ and $w(3,2)\sspeq \overline w \sspeq -\frac{1}{2}(1\sspp \sqrt{3}\,i)$ one finds \pn
$[V^{-1}_3(x)]_{1,2}\sspeq  w/(3\,x)$, because $DN(3,2,x) = x^2\,(1\sspm \overline w)\,(w\sspm \overline w)\sspeq x^2\,\frac{1}{2}\,(3\sspp i\,\sqrt{3})\,i\,\sqrt{3}  \sspeq  3\,w\,x^2$, and from $P(3,2) = [1,\, w]$, and $Ch(3,1,2)\sspeq [[1],\, [w]]$ one obtains $N(m,l,j,x)\sspeq (-1)^1\,x\,(1\sspp w) = x\,\overline w$. Indeed,  $[V^{-1}_3(x)]_{1,2} \sspeq \overline w/(3\,w)\sspeq w/(3\,x)$, due to $\overline w^2 \sspeq w$.\psn
\Beq
V_3^{-1}(x)\sspeq \frac{1}{3}\,\LTmatd{1}{1}{1}{1/x}{\overline w/x}{w/x}{1/x^2}{w/x^2}{\overline{w}/x^2}\,.
\Eeq
Therefore the standard trisection of $G(x)$ is
\Beqarray
G(3,0,x) &\sspeq & \frac{1}{3}\,\left (G(x^{1/3}) \sspp G(w\,x^{1/3}) \sspp  G(\overline w\,x^{1/3})\right), \\
G(3,1,x) &\sspeq & \frac{1}{3\,x}\,\left(G(x^{1/3}) \sspp \overline w\,G(w\,x^{1/3}) \sspp w\, G(\overline w\,x^{1/3})\right), \\
G(3,2,x) &\sspeq & \frac{1}{3\,x^2}\,\left(G(x^{1/3}) \sspp w\, G(w\,x^{1/3}) \sspp  \overline w\, G(\overline w\,x^{1/3})\right)\,.
\Eeqarray
This should then be simplified for given $G(x)$, by finding the rational function $P(x)/Q(x)$ which can become tedious in the general m-section case (the computer will help). \pn
The topic of this paper is to give for the general linear three-term recurrence relation the coefficients of these polynomials $P$ and $Q$ in terms of well known polynomials which are functions of the signature of this recurrence.    
\pbn 
\section {\bf General linear three term recurrence} 
This section is a review of basic formulas for the considered recurrence relation.\pn
The sequence $\{H(p,q;r,s;n)\}_{n=0}^{\infty}$ satisfies the following linear three- term (also called second order) recurrence relation of signature $(r,\,s)$, with integer numbers $r$ and $s$, both non-vanishing, and initial conditions (seeds or inputs) $(p,\,q)$, with integer numbers $p$ and $q$. Only integer sequences are considered. In the following these domains for $p,\,q,\,r,\,s$ will not be repeated in the formulas.\pn
The letter $H$ is used because this sequence has been studied by {\sl A. F. Horadam} in many publications. See \eg \cite{Horadam1}, \cite{Horadam2},\cite{Horadam3}, and also \cite{MathWorld}. 
\Beqarray
 H(p,q;r,s;\,n) &\sspeq& r\,H(p,q;r,s;\,n-1) \sspp  s\,H(p,q;r,s;\,n-2),\ \text{for}\ n \sspgeq  2,\ {\text{and}}  \\ 
 H(p,q;r,s;\,0) &\sspeq& p,\,\  H(p,q;r,s;\,1) \sspeq q.
\Eeqarray
It is sufficient to consider the seeds $(p,\,q) \sspeq (0,\,1)$, naming the sequence $\{H01(r,s;\,n)\}_{n=0}^{\infty}$, because 
\Beq
H(p,q;r,s;\,n) \sspeq q\,H01(r,s;\,n) \sspp  p\,s\,H01(r,s;\,n-1).
\Eeq
Also $H01(r,s;-1) \sspeq 1/s$ and $H01(r,s;-2) \sspeq -r/s^2$ will be used. \pn 
One can also extend this sequence to all negative integer $n$, by
\Beq
H01(r,s;\,n) \sspeq -(-s)^n\, H01(r,s;\,-n),
\Eeq
which implies the result for negative indices for sequence $H$.
\psn
The {\sl Binet} - {\sl de Moivre} formula is
\Beq
H01(r,s;n)\sspeq \frac{\lambda(r, s)^n \sspm (-s/\lambda(r, s))^n} {\lambda(r, s) \sspm (-s/\lambda(r, s))},\ \text{where}\ \lambda(r,\,s)\sspeq \frac{1}{2}\,\left (r\sspm \sqrt{r^2\sspp 4\,s}\,\right)\,.
\Eeq
The transfer matrix, also called $\bf Q$ matrix, for the $(r\, s)$ recurrence relation is
\Beq
{\bf Q}\sspeq \LTmatz{r}{s}{1}{0}\,.
\Eeq
The powers of this $2\times2$ matrix with trace $Tr\,{\bf Q}\sspeq r$ and determinant $Det\,{\bf Q}\sspeq -s \sspneq 0$, can be found with the help of the {\sl Cayley}-{\sl Hamilton} Theorem in terms of Chebyshev $S$-polynomials by
\Beq
{\bf Q}^n(r, s) \sspeq (\sqrt{-s})^n\,\left[ S\left(n,\,\frac{r}{\sqrt{-s}}\right)\,{\bf 1}  \sspp S\left(n-1,\,\frac{r}{\sqrt{-s}}\right)\,\frac{1}{\sqrt{-s}}\,({\bf Q}(r,s)\sspm  r\,{\bf 1})\right]\,. 
\Eeq
For the {\sl Chebyshev} S-polynomials see OEIS \cite{OEIS} \seqnum{A049310} for their coefficients, their properties, and references, \eg \cite{ASt}, \cite{Rivlin}. OEIS $A$-number links will henceforth be used without citation.
\Beq
S(n,\,x) \sspdef H(1,x;x,-1;\,n),\ {\text {for}}\ n\sspgeq 0. 
\Eeq

For negative $n$ one finds $S(-1,\, x) \sspeq 0$, and $S(n,\,x) \sspeq -S(-n-2,\,x)$, for $n \sspleq -2$. 
\psn
This produces the matrix
\Beq
{\bf Q}^n(r,\,s) \sspeq (\sqrt{-s})^n\,\LTmatz{S\left(n,\,\frac{r}{\sqrt{-s}}\right)}{\frac{s}{\sqrt{-s}}\,S\left(n-1,\, \frac{r}{\sqrt{-s}}\right)}{\frac{1}{\sqrt{-s}}\,S\left(n-1,\, \frac{r}{\sqrt{-s}}\right) }{-S\left(n-2,\, \frac{r}{\sqrt{-s}}\right)}\,.
\Eeq
The (generalized) {\sl Chebyshev} $T$-polynomials are defined from the trace as
\Beq
T\left(n,\, \frac{r}{2\,\sqrt{-s}}\right)\sspdef \frac{1}{2}\,Tr\,{\bf Q}^n(r,\,s) \sspeq \frac{1}{2}\,\left( S\left(n,\, \frac{r}{\sqrt{-s}}\right) -  S\left(n-2,\, \frac{r}{\sqrt{-s}}\right) \right)\ .
\Eeq
For $(r,s)\sspeq (x,-1)$ these are the usual {\sl Chebyshev} $T$-polynomials: $T(n,\,x/2)\sspeq \frac{1}{2}\,(S(n, x) \sspm S(n-2,\,x))$. Later $R(n,\, x) \sspeq T(n,\,x/2)/2$ will be used.
\pn
Because $Det\,{\bf Q}(s,\,r)\sspeq -s$ one has $Det\,{\bf Q}^n(s,\,r)\sspeq (-s)^n$, by the product theorem for determinants, and this leads to the {\sl Cassini}-{\sl Simson}
identity in $n$ and $(r,s)$ (with $n \to n+1$)
\Beq
S(n,\, y)^2\sspm S(n-1,\,y)\,S(n+1,\,y)\sspeq 1, 
\Eeq
where $r$ and $s$ only enter \via\, $y\sspeq y(r,\, s) \sspeq \frac{r}{\sqrt{-s}}$. 
\pbn
A further reduction of the $H01$ sequence, important for the main part of this paper, is possible in terms of the usual {\sl Chebyshev} $S$-polynomials by
\Beq
H01(r,s;n) \sspeq (\sqrt{-s}\,)^{n-1} \,S\left(n-1,\,\frac{r}{\sqrt{-s}}\right)\,.
\Eeq
This follows from comparing the recurrence and the seeds. 
\pbn
The ordinary generating functions (\ogf\llap) of $\{H01(r,s;n\}_{n=0}^{\infty}$ is
\Beq
GH01(r,s;x)\sspeq \frac{x}{1\sspm r\,x \sspm s\,x^2}\,.
\Eeq
The \ogf of $\{H(p,q;r,s;n\}_{n=0}^{\infty}$ in terms of $GH01(r,s;x)$ is
\Beqarray
GH(p,q;r,s;x)&\sspeq& p + (q\sspp p\,s\,x)\,GH01(r,s;x),\\
 &\sspeq& \frac{p \sspm (p\,r \sspm q)\,x}{1\sspm r\,x \sspm s\,x^2}\,.
\Eeqarray
The \ogf of $\{S(n,\,y\}_{n=0}^{\infty}$ is 
\Beq
GS(x)\sspeq \frac{1}{1\sspm y\,x \sspp x^2}\,.
\Eeq
\pbn
\section {\bf Conjecture for improved formulas for the m-section of the linear three -term recurrence sequences} 
\hskip .5cm This section contains conjectures for simplified formulas for the $m$-section or the special sequences $H$, $H01$ and $S$. In the next section these conjectures will be proved.\psn
One of the authors ({\sl G.\, D.}) heuristically found a formula for the sequence $\{H(p,q;r,s;\, m\,n)\}_{n=0}^\infty$ , for $m\sspgeq 0$, that identifies it as an $H$ sequence with different input $(p,\, q^\prime)$ and signature $(r^\prime,\, s^\prime)$. See his comment in \seqnum{A034807} where $p,\ q,\,r,\, s$ are denoted as $a,\,b,\,c,\,d$, respectively.
\pn
The second author generalized this conjecture to the $m$-section of the sequence $H$ and their \ogf\,\llap s. He also proved a conjecture for the sequence $H01$ which implies the one for $H$. In the next section the proof will be given for the conjecture for the $m$ section of {\sl Chebyshev} $S$-polynomials and the \ogf\,\llap s., that will lead to the other two conjectures.\
\psn
{\bf Conjecture for $H$}
\Beqarray
&& \rm{For}\ n\sspgeq 0,\ m\sspgeq (1),2\ \rm{and}\ l\sspin\{\range{0}{m-1}\}\nonumber \\   
&& H(p,q;r,s;\,m\,n+l) \sspeq H(H(p,q;r,s;\,l),H(p,q;r,s;\,m+l);SUM(r,s;\,m),-(-s)^m;\,n)\,.
\Eeqarray
with
\Beq
SUM(r,s;\,m)\sspeq r^m\,\sum_{k=0}^{\floor{\frac{m}{2}}}\,\seqnum{A034807}(m,\, k)\,(s/r^2)^k\,.
\Eeq
Therefore $SUM(r,s;\,m)$ is the polynomial $P(m,\,x)$ of row $m$ of this irregular triangle evaluated at $x\sspeq s/r^2$ and scaled by $r^m$.
\psn
Note that the symmetry between $n$ and $m$ for the \lhs is not obvious for the \rhs, but true because the later proof can be done with interchanged $n$ and $m$.  This symmetry holds for all versions of the conjecture given later.
\psn
The recurrence relation for the triangle $T$\sspeq \seqnum{A034807} (given there by {\sl Michael Somos}, given here without proof) is 
\Beqarray
&& T(n,\,k) \sspeq T(n-1,\, k) + T(n-2,\, k-1), \ \rm{for}\  n\sspgeq 2,\ \rm{and} \nonumber \\
&& T(0,\,k) \sspeq 2,\ \text{for}\ k\sspeq 0,\ \text{otherwise}\ 0,\nonumber \\
&& T(1,\,k) \sspeq 1,\ \text{for}\ k\sspeq 0,\ \text{otherwise}\ 0. 
\Eeqarray
The explicit form (given in \seqnum{A034807} by {\sl Alexander Elkins}, here also given without proof) is
\Beq
T(n,\,k) \sspeq \frac{n\,(n-1-k)!}{k!\,(n\sspm 2\,k)!},\ \rm{for}\  n \sspgeq 1, \, k\sspeq 0,\,1,\,...,\, \floor{n/2}\,\ \rm{and}\  T(0,\, 0) \sspeq 2\,.  
\Eeq
The \ogf for the row polynomials $\{P(n,\,y)\}$ of $T$ (the \ogf of the triangle) (given there by {\sl Vladeta Jovovic}, here given also without proof) is
\Beq
\frac{2 \sspm x}{1 \sspm x \sspm y\,x^2}\,.
\Eeq
{\bf Lemma 1}
\Beqarray
{\bf i)}\ SUM(r,s;\,m)&\sspeq& H(2,r;r,s;\,m), \\
{\bf ii)}\phantom{\ SUM(r,s;\,m)}  &\sspeq& s\,H01(r,s;\,m-1) \sspp H01(r,s;\,m+1),\\
{\bf iii)}\phantom{\ SUM(r,s;\,m)} &\sspeq& (\sqrt{-s}\,)^m\,R\left(m,\,\frac{r}{\sqrt{-s}}\right), \ \rm{with}\ R(n,\, x)\sspdef S(n,\, x) \sspm S(n-2,\,x)\ .
\Eeqarray
The polynomials $R$ are the monic {\sl Chebyshev} $T$-polynomials. See \seqnum{A127672}
 for their coefficients and properties.
\psn
{\bf Proof}
\psn
{\bf i)} This follows from the definition of $SUM$  in \Eq{31} and the recurrence of the irregular triangle $T$ given in \Eq{32}, leading to the signature $(r,\,s)$, and the inputs $SUM(r,s;\, 0) \sspeq 2$ and $SUM(r,s;\, 1) \sspeq r$. \psn  
{\bf ii)} Results after replacing the sequence $H$ by $H01$ according to \Eq{16}.
\psn
{\bf ii)} Uses the replacement of the sequence $H01$ by the {\sl Chebyshev} $S$-polynomials , \Eq{25}, evaluated at $x\sspeq r/\sqrt{-s}$. \hskip 15.7cm $\square$
\psn 
{\bf Example 2: Fibonacci trisection}\psn
$ F(n)\sspeq H(0,1;1,1;\, n)\sspeq$\seqnum{A000045}$(n)$, for $n\sspgeq 0$. The first part of the trisection ($m\sspeq 3,\,l\sspeq 0$) is $F(3\,n)\sspeq$\seqnum{A014445}$(n)\sspeq \{0,\,2,\, 8,\, 34,\, 144,...\}$. The conjecture leads to \pn
$F(3\,n) \sspeq H(0,\, F(3);\, -i\,R(3,\,1/\sqrt{-1}\,),-(-1)^3;\,n)$, where $ R(3,\,-i) \sspeq 4\,i$. $F(3\,n)\sspeq H(0,2;4,1;\,n)$. 
\pn
In terms of $H01$ this becomes $F(3\,n)\sspeq 2\,H01(4,1;\,n)$.\pn
In terms of $S$ one finally finds $F(3\,n)\sspeq 2\,i^{n-1}\,S(n-1,-4\,i)$.\psn
The other parts of the trisection $F(3\,n\sspp 1)\sspeq$\seqnum{A033887}$(n)$ and $F(3\,n\sspp 2)\sspeq$\seqnum{A015448}$(n+1)$, for $n\sspgeq 0$,follow similarly, and the results are
\Beqarray
F(3\,n)&\sspeq& H(0,2;4,1;\,n) \sspeq 2\,H01(4,1;\,n) \sspeq 2\,i^{n-1}\,S(n-1,-4\,i),   \\
F(3\,n\sspp 1)&\sspeq& H(1,3;4,1;\,n) \sspeq H01(4,1;\,n+1) \sspm H01(4,1;\,n)\\ &&  \hskip 2.5cm \sspeq -i^n\,(S(n,\,-4\,i) \sspp i\,S(n-1,\,-4\,i)),\nonumber \\
F(3\,n\sspp 2)&\sspeq& H(1,5;4,1;\,n)\sspeq H01(4,1;\,n+1)\sspp H01(4,1;\,n)\\
&&  \hskip 2.5cm \sspeq i^n\,(S(n,\,-4\,i) \sspm i\,S(n-1,\,-4\,i))\,.\nonumber
\Eeqarray
\psn
The  above conjecture for $H$ is equivalent to the one for its \ogf 
\Beq
GHml(p,q;r,s;m,l;\, x)\sspdef \sum_{n=0}^\infty\,H(p,q;r,s;\,m\,n\sspp l)\,x^n\,.
\Eeq
\psn
{\bf Conjecture for $\bf GHml$}\psn
\Beq
GHml(p,q;r,s;m,l;\,x) \sspeq \frac{H(p,q;r,s;\,l) \sspm \left(H(p,q;r,s;\,l)\,SUM(r,s;\,m) \sspm H(p,q;r,s;\,m+l)\right)\,x}{1\sspm SUM(r,s;\,m)\,x\sspp (-s)^m\,x^2}\,.
\Eeq
{\bf Proof:} This equivalence of conjectures is clear from the \ogf of the sequence $H$ given in \Eq{28}. One has just to insert the conjectured values for the inputs and signature from \Eq{30}\hskip 4.5cm $\square$
\psn
{\bf Example 3: {\bf O.g.f.} Fibonacci trisection}\psn
For $ m \sspeq 3$, $(p,\,q)\sspeq (0,\,1)$ and $(r,\,s)\sspeq (1,\,1)$ the denominator of $GF3l(x) \sspdef \sum_{n=0}^\infty\,F(3\,n+l)\,x^n$ is $1\sspm (-i\,R(3,\,1/i))\,x \sspp (-1)^3\,x^2\sspeq  1 \sspm 4\,x\sspm x^2$, for $l\sspeq 0,\,1$ and $2$. The numerators are $F(3)\,x \sspeq 2\,x$, $1 \sspm (1\cdot 4\sspm 3)\,x\sspeq 1\sspm x$, and $1 \sspm (1\cdot 4\sspm 5)\,x\sspeq 1\sspp x$, for these $l$ values, respectively.
\psn
The conjecture for the $m$-section of $H$ implies the one for $H01$, and the corresponding \ogf \llap s, are obtained setting $(p,\,q)\sspeq (0,\,1)$, and then rewriting $H$ in terms of $H01$ using \Eq{16} with the new parameters.
In Example 2 this second step has been used for $m\sspeq 3$ and $(r,\,s) \sspeq (1,\,1)$.
\psn
{\bf Conjecture for $\bf H01$}
\Beqarray
&& H01(r,s;\,m\,n\sspp l)  \sspeq  q^\prime H01(r^\prime,s^\prime;\,n) \sspp p^\prime\,s^\prime\,H0(r^\prime,s^\prime;\,n-1)\,,\ \rm {where} \nonumber \\
&& p^\prime \sspeq p^\prime(r,s;\,l) \sspeq H01(r,s;\,l),\, \ q^\prime \sspeq p^\prime(r,s;\,\,m+l) \sspeq H01(r,s;\,m+l),\, \nonumber \\ 
&& r^\prime\sspeq r^\prime(r,s;\,m)\sspeq (\sqrt{-s}\,)^m\,R(m,\,r/\sqrt{-s}\,),\ \ s^\prime \sspeq s^\prime(s,\,m)\sspeq  -(-s)^m. 
\Eeqarray
\psn
The part $l\sspeq 0$ simplifies to
\Beq
H01(r,s;\,m\,n)\sspeq H01(r,s;\,m)\,H01((\sqrt{-s}\,)^m\,R(m,\,r/\sqrt{-s}\,),\,-(-s)^m;\,n)\,.
\Eeq
\psn
{\bf Conjecture for $\bf GH01ml$}
\psn
The conjecture for the \ogf $GH01ml(r,s;m,l;\,x)\sspdef \sum_{n=0}^\infty H01(r,s;m\,n+l)\,x^n$ is obtained from $GHml(x)$ in \Eq{42}, and is given with $ y\sspeq y(r,\,s)/\sqrt{-s}$ as   
\Beq
GH01ml(r,s;m,l;\,x)\sspeq \frac{H01(r,s;\,l)\sspm ((\sqrt{-s}\,)^m\,H01(r,s;\,l)\,R(m,\,y) \sspm H01(r,s;\,m\sspp l))\,x}{1\sspm (\sqrt{-s}\,)^m\,R(m,\, y)\,x \sspp (-s)^m\,x^2 }\,. 
\Eeq
\psn
Because the sequences $H$ and $H01$ are determined by the {\sl Chebyshev} polynomials $\{S(n,\, y\sspeq r/\sqrt{-s}\,)\}$ the conjecture for $S(m\,n\sspp l,\,r/\sqrt{-s}\,)$ is fundamental.
\psn
\vfill\eject
\noindent
{\bf Conjecture for $\bf S$}
\psn
For $n \sspgeq 0,\,m\sspgeq (1),2$ and $l\sspin\{\range{0}{m-1}\}$:  
\Beqarray 
S(m\,n\sspp l,\,y) &\sspeq& c(s,\,m)^{n-1}\,\left\{S(m\sspp l,\,y)\, S(n-1,\,c(s,\,m)\,R(m,\,y)) \right. \,\nonumber \\ 
&& \hskip 2cm \left. \sspm c(s,\,m)\,S(l,\,y)\, S(n-2,\,c(s,\,m)\,R(m,\,y)) \right\},\\ 
&& \rm{with}\ y\sspeq r/\sqrt{-s},\, c(s,\,m)\sspdef \frac{(\sqrt{-s}\,)^m}{\sqrt{(-s)^m}},\,S(-2,y)\sspeq -1,\ \rm{and}\ S(-1,y)\sspeq 0\,. \nonumber
\Eeqarray
The part $l\sspeq 0$ simplifies, using $c(s,\,m)^2\sspeq 1$, the recurrence relation of $S$ and then the definition of $R$, to
\Beqarray
&&S(m\,n,\,r/\sqrt{-s})\sspeq (c(s,\,m))^n\,{\bf \cdot}\nonumber \\
&&{\bf \cdot} \{S(n,\,c(s,\,m)\,R(m,\,r/\sqrt{-s})) + c(s,\,m)\,S(m-2,\,r/\sqrt{-s}\,)\, S(n-1,\,c(s,\,m)\,R(m,\,r/\sqrt{-s}\,))\}.
\Eeqarray
\psn
The following proof that this conjecture is equivalent to the conjecture for $H01$ uses 
 $y\sspeq r/\sqrt{-s}$ and \Eq{25}.
 \psn
{\bf Proof of the equivalence between the conjectures $\bf H01$ and $\bf S$}
\psn
With $ y\sspeq r/\sqrt{-s}$ and \Eq{25}
$S(m\,n\sspp l,\,y) \sspeq (1/\sqrt{-s}\,)^{m\,n\sspp l})\,H01(r,s;m\,n \sspp l \sspp 1)$. With the conjecture for $H01$ from above this becomes in terms of $S$, again using \Eq{25},\pn
\Beqarray
(\sqrt{-s}\,)^{m\,n\sspp l}\,S(m\,n\sspp l,\,y) &\sspeq& \widehat{q^\prime}\,(\sqrt{-s'}\,)^{n-1}\,S(n-1,r^\prime/\sqrt{-s^\prime}\,) \nonumber \\
&&\sspp \widehat {p^\prime}\,s^\prime\,(\sqrt{-s^\prime}\,)^{n-2}\, S(n-2,r^\prime/\sqrt{-s^\prime}\,)\,,\nonumber
\Eeqarray
with $r^\prime,\,s^\prime,\,p^\prime$ and $q^\prime$ from \Eq{43},\, and $q^\prime$ and $p^\prime$ are written in terms of $S$ as $\widehat{q^\prime} \sspeq \sqrt{-s}\,^{m+l}\,S(m+l,\, y)$ and $\widehat{p^\prime}\sspeq \sqrt{-s}\,^{l}\,S(l,\, y) $. Also $r^\prime/\sqrt{-s^\prime}\sspeq c(s,m)\,R(m,\,y)$.\pn
Dividing both sides by $(\sqrt{-s}\,)\,^{m\sspp l}\, (\sqrt{(-s)^m}\,)\,^{n-1}$ produces
\Beqarray
c(s,m)^{n-1}\,S(m\,n\sspp l,\,y) &\sspeq& \{S(m\sspp l,\,y)\,S(n-1 \, c(s,m)\, R(m\,y)) \nonumber \\
&&\sspm (1/c(s,m))\,S(l,\,y)\, S(n-2,\,c(s,m)\, R(m,\,y)) \}\ . 
\Eeqarray
Because $c(s, m)^2\sspeq 1$ one replaces $1/c(s,\,m)$ by $c(s,\,m)$, giving the final result. \hskip 4cm  $\square$
\psn
Note that $c(s,\,m)$ has only values from $\{+1,\,-1\}$. $c(-1,\,m) \sspeq 1$, for $m\sspgeq 1 $, 
and $\{c(1,\,m)\}_{m\sspgeq 1} \sspeq$ repeat $\{1,\,-1,\,-1,\,1\} \sspeq$ \seqnum{A087960} with offset $1$.
\psn
{\bf Example 4: Trisection of Chebyshev S-polynomials} \psn
 m\sspeq 3$, $r\sspeq y,\,s\sspeq -1. Note that $y$ is now an indeterminate.\pn
${\bf l\sspeq 0}:\  S(m\,n,\, y) \sspeq S(n,\, R(3,\, y)) \sspp y\,S(n-1,\, R(3,\, y))$,
with $R(3,\,y)\sspeq y\,(y^2\sspm 3)$.\pn
${\bf l\sspeq 1}:\ S(m\,n\sspp 1,\, y) \sspeq S(4,\,y)\,S(n-1,\,R(3,\, y)) - y\,S(n-2, R(3,\, y)))$, with $S(4,\,y)\sspeq 1\sspm 3\,y^2\sspp y^4$.\pn 
${\bf l\sspeq 2}:\ S(m\,n\sspp 2,\,y) \sspeq S(2,\,y)\,(R(3,\,y)\,S(n-1,\, R(3,\,y)) \sspm S(n-2,\, R(3,\,y))\sspeq S(2,\,y)\,S(n,\, R(3,\,y))$,\pn 
because $S(5,\,y)\sspeq S(2,\,y)\,R(3,\,y)$, and $S(2,\,y)\sspeq y^2\sspm 1$.
\pbn
The conjecture for the \ogf $GSml(r,\,s;m,l;\,x)\sspdef \sum_{n=0}^\infty S(m\,n\sspp l,\,y\sspeq r/\sqrt{-s}\,)\,x^n$ is obtained from the one for $GHml(x)$ given above.
\psn
{\bf Conjecture for $\bf GSml$}\psn
With $y\sspeq \frac{r}{\sqrt{-s}}:$
\Beqarray
GSml(r,s;m,l;\,x)&\sspeq& \frac{1}{(\sqrt{-s}\,)^l}\,GH01ml\left(r,s;m,l\sspp 1;\,\frac{x}{(\sqrt{-s}\,)^m}\right), \nonumber \\
&\sspeq& \frac{S(l,\,y) \sspm (S(l,\,y)\,R(m,\,y)
\sspm S(m\sspp l,\,y))\,x}{1\sspm R(m,\,y)\,x \sspp x^2}\,.
\Eeqarray
Note that the advantage of working with the \ogf\llap s instead of the sequences is that the $(r,\,s)$ dependence appears only in $y$ (not like in \Eq{46} also in $c(s,\,m)$). 
\psn
{\bf Exercise}\psn
In order to appreciate these formulas on should compare them with the standard computation according to {\it Section 2}. Done either by hand or by computer the result will not be expressed in terms of {\sl Chebyshev} polynomials.
\psn
{\bf Proof of the equivalence between $\bf GSml$ and $\bf GH01ml$}\psn
This uses the relation between $S(n,\,y)$ and $H01(r,s;\, n+1)$ obtained from \Eq{25}, for $n\sspto m\,n\sspp l$. This leads to the relation between the \ogf\llap s. Then in \Eq{45} the $H01$ sequences are rewritten in terms of $S$, with $y\sspeq r/\sqrt{-s}$.\hskip 14.2cm $\square$
\psn
{\bf Example 5:\ O.g.f.s for the trisection of Chebyshev $\bf S$ polynomials}\psn
 $m\sspeq 3$,\, $r\sspeq y,\,s\sspeq -1$. Note that $y$ is now an indeterminate.\pn
${\bf  l\sspeq 0}:\ GS30(y,\,x) \sspeq (1 - (R(3,\,y) \sspm S(3,\,y))\,x)/(1\sspm R(3,\,y)\,x\sspp x^2)$, With $R(3,\, y)$ from above in Example $4$, and $S(3,\,y)\sspeq y\,R(2,\,x)\sspeq y\,(y^2\sspm 2)$ one obtains $R(3,\,y) \sspm S(3,\,y) = -y$, hence $GS30(y,\,x)\sspeq (1\sspp y\,x)/(1\sspm y\,(y^2\sspm 3)\,x\sspp x^2)$.
\pn
${\bf l\sspeq 1}$:\ In the numerator appears  $y\,R(3,\,y)\sspm S(4,\,y)\sspeq -1$. Hence\pn $GS31(y,\,x) \sspeq (y\sspp x)/(1\sspm y\,(y^2\sspm 3)\,x\sspp x^2)$.\pn
${\bf l\sspeq 2}$: In the numerator appears  $S(2,\,y)\,R(3,\,y)\sspm S(5,\,y)\sspeq 0$ (see Example $4$). Hence\pn
$GS32(y,\,x) \sspeq (y^2\sspm 1)/(1\sspm y\,(y^2\sspm 3)\,x\sspp x^2)$.\pn
\psn
\section {\bf Proof of the conjectures}\psn
The proof is given for the conjectured \ogf\llap s, equivalent to the conjectures for the corresponding sequences.\pn
Here the proof for the conjecture of the \ogf of the sequence $S$, \ie $GSml$ of \Eq{49}, is given  which is equivalent to the \ogf of sequence $H01$, \ie $GH01ml$ of \Eq{45}.\pn
The conjecture for the \ogf of the sequence $H$, \ie  $GHml$ of \Eq{42}, follows from the 
conjecture of $GSml$ by 
\Beqarray
GHml(p,q;r,s;m,l;\,x)&\sspeq& q\,(\sqrt{-s}\,)^{l-1}\,GSml\left(r,s;m,l-1;\,(\sqrt{-s}\,)^m\,x\right) \nonumber \\
&& \sspp p\,s\,(\sqrt{-s}\,)^{l-2}\,GSml\left(r,s;m,l-2;\,(\sqrt{-s}\,)^m\,x\right)\,.
\Eeqarray
Note that for $m\sspgeq 2$ and $l\sspeq 0$ and $1$ the sequences $H$, $H01$ and $S$ appear also with negative indices $n\sspeq -1$ and $-2$, namely $H01(r,s; -1)\sspeq 1/s$, $S(-2,\,x)\sspeq -1$ and $ S(-1,\,x)\sspeq 0$.\psn
This $GHml$ formula coincides with the original one of \Eq{42}, after sequence $H$ is replaced by sequence $H01$, and then by sequence $S$.
\psn
{\bf Theorem:} {\it The conjecture for the \ogf $GSml$ of $\{S(m\,n+l,\,r/\sqrt{-s}\,)\}_{n=0}^\infty$ is true. }
\psn
{\bf Proof:}\psn
One proves that the \ogf $GS(y,\,x)\sspeq 1/(1\sspm y\,x\sspp x^2)$ for the {\sl Chebyshev} polynomials $\{S(n,\,y)\}_{n=0}^\infty$, with $y \sspeq r/\sqrt{-s}$, satisfies the $m-$section formula according to \Eq{1} in terms of the conjectured part $l$ \pn \ogf\llap s $GSml$ from \Eq{49}. This can be rewritten, by bringing the identical ($l$-independent) denominators of $GSml$ to the left hind side, and the denominator of $GS(y,\,x)$ to the right-hand side as
\Beqarray
LHS(m,l;\, y,\,x)&\sspdef& 1\sspm R(m,\,y)\,x^m\sspp x^{2\,m}, \nonumber\\
RHS(m,l;\, y,\,x)&\sspdef& (1\sspm y\,x\sspp x^2)\,\sum_{l=0}^{m-1}\,x^l\, N(m,l;\,y,\,x^m), \nonumber\\
\rm{with}\ && N(m,l;\,y,\,x^m)\sspeq S(l,\,y)\sspm (S(l,\,y)\,R(m,\,y)\sspm S(m+l,\,y))\,x^m\,.
\Eeqarray
 Remember that by working with \ogf\llap s instead of sequences the $(r,\,s)$ dependence appears only $y$. Therefore the proof will be given for the indeterminate y.
 \psn
 Because the {\sl Vandermonde} matrix has an inverse (see \Eq{4}) the proof will automatically hold also for $GSml$ in terms of $GS$ like in \Eq{9}. 
 \psn
 All powers of $x$ will be compared on both sides in order to prove that $LHS \sspeq RHS$.
 \pn
 In $RHS$ all powers $x^0,\,x^1,\, ...,\, x^m,\,...,\,x^{2\,m},\,x^{2\,m+1}$ appear. \
 In $LHS$ only $x^0,\,x^m,\,x^{2\,m}$ are present.
 \pn
 It will turn out that the proof for the two highest powers $x^{2\,m+1}$ and $x^{2\,m}$ differs from the one for the other powers. Usually the recurrence of the {\sl Chebyshev} $S$ polynomials will show directly that $RHS\sspm LHS \sspeq 0$ but for the two highest powers one has to use results from the bisection of these polynomials.
 \pn
 For the other powers the contribution of the $R$ terms in the numerator $N$ of $GSml$ will be considered separately from the remainder (pure $S$ terms). For the powers $x^m$ to $x^{2\,m-1}$ it will turn out that the  $R$ terms are multiplied by factors which vanish because of the recurrence of the $S$-polynomials (the structure of $R$ will thus be irrelevant).
 \psn 
 One starts with the two highest powers.\pn 
 For $x^{2\,m+1}$ only $RHS$ is present, namely $-(S(m-1,\,y)\,R(m,\,y) \sspm S(2\,m\sspm 1))$. It vanishes if 
 \Beq
 S(2\,m\sspm 1,\,y)\sspeq S(m-1,\,y)\,R(m,\,y)\,.
 \Eeq
 For $x^{2\,m}$ the $RHS$ becomes (-y)\,\{-(S(m-1,\,y)\,R(m,\,y) \sspm S(2\,m\sspm 1,\,y))\}\sspp (+1)\,\{-(S(m-2,\,y)\,R(m,\,y) \sspm \pn
 S(2\,m\sspm 2,\,y))\}, and $RHS \sspeq 1$. The first term vanishes if the $x^{2\,m+1}$ power contribution vanishes, and then for this $x$-power $RHS\sspm LHS\sspeq 0$ if 
 \Beq
 S(2\,(m-1),\,y)\sspeq 1 \sspp S(m-2,\,y)\,R(m,\,y)\,.
 \Eeq
 {\it {\bf Lemma 2}: Eqs. (52) and (53) are satisfied for all $m\sspgeq 0$\,.} 
 \psn
 {\bf Proof:}\pn
 These two equations are found in \cite{MOS}, written for Chebyshev $T$ and $U$-polynomials.
 \pn
 Using a proof of the standard bisection will not help here. The proof is done by induction on $m$ on both equations simultaneously, employing the {\sl Cassini-Simson} identity from \Eq{24}.\pn
 For $ m\sspeq 0 $ the first equation is fulfilled because $S(-1, \,y) \sspeq 0$, and the second one because  $S(-2,\, y)\sspeq -1$ and $R(0,\,y)\sspeq 2$.
 \pn
 Assume that both equations hold for $m'\sspeq 1,\,2,\, ...,\,m$. First, \Eq{53} will be proved for $m\sspto m+1$. Multiplying \Eq{52} by $y$ and subtracting \Eq{53} yields, after using recurrence relations, 
 \Beq
 S(2\,m,\,y)\sspeq S(m,\,y)\,R(m,\,y) \sspm 1\,.
 \Eeq
 For \Eq{53} one wants to prove $S(2\,m,\,y)\sspeq 1\sspp S(m-1,\,y)\,R(m+1,\,y)$, \ie $ 0 \sspeq 1\sspp S(m-1,\,y)\,R(m+1,\,y) - (S(m,\,y)\,R(m,\,y) \sspm 1)$, \ie $0\sspeq 2\sspm S(m,\,y)\,R(m,\,y)\sspp  S(m-1,\,y)\,R(m+1,\,y) $. Replacing $R$ in terms of $S$, using  $ S(m\sspm 1,\,y)\, S(m\sspp 1,\,y) \sspeq S(m,\,y)^2\sspm 1$ ({\sl Cassini-Simson})  leads to a cancellation of $S(m,\,y)^2$, leaving  $0 \sspeq 1 \sspm S(m\sspm 1,\,y)^2  \sspp S(m,\,y)\,S(m\sspm 2,\,y)$, which is again a {\sl Cassini-Simson} identity. Thus \Eq{53} is proved.
 \psn
 For \Eq{52} one wants to prove $S(2\,m\sspp 1,\,y)\sspeq S(m,\,y)\,R(m\sspp 1,\,y)$. By recurrence $S(2\,m\sspp 1,\,y)\sspeq y\,S(2\,m,\,y) \sspm S(2\,m\sspm 1,\,y)$ which becomes, with the induction assumptions \Eq{54} and \Eq{52} (for $m$), $S(2\,m\sspp 1,\,y)\sspeq  -y \sspp y\,S(m,\,y)\,R(m,\,y)\sspm S(m\sspm 1,\,y)\,R(m,\,y)$. This is by recurrence $S(2\,m\sspp 1,\,y)\sspeq -y \sspp S(m\sspp 1,\,y)\,R(m,\,y)$. One wants now to prove $S(2\,m\sspp 1,\,y) \sspeq S(m,\,y)\,R(m+1,\,y) \sspeq -y \sspp S(m\sspp 1,\,y)\,R(m,\,y)$. Replacing $R$ in terms of $S$  gives, after cancellation of  $S(m,\,y)\, S(m\sspp 1,\,y)$, $S(m\sspp 1,\,y)\,S(m\sspm 2,\,y)\sspm S(m,\,y)\,S(m\sspm 1,\,y)\sspeq -y$. Replacing  $S(m\sspm 2,\,y)$ by $y\,S(m\sspm 1,\,y)\sspm S(m,\,y)$, and then $S(m\sspp 1,\,y)\,S(m\sspm 1,\,y) \sspeq -1\sspp S(m,\,y)^2$ ({\sl Cassini-Simson)} leads to $ 0\sspeq S(m,\,y)\,(y\, S(m,\,y) \sspm S(m\sspm 1,\,y))\sspm S(m\sspp 1,\,y)\,S(m,\,y)$, which holds because of the recurrence relation. \hskip 3cm $\square$ 
 \psn
Consider the $R$ term contributions in the numerator $N$ together with the pre-factor with powers $x^i$ for $i\sspin \{0,\,1\,,2\}$. The powers are $x^{i+l+m}$, for $m\sspgeq 2$ and $l\sspin\{\range{0}{m-1}\}$, but only for $e\sspdef i\sspp l\sspp m <= 2\,m\sspm 1$, because the powers $x^m$ and $x^{2\,m+1}$ have just been treated separately. $R$ terms appear only for the exponents $e\sspeq 2\,m\sspm 1,\,...,\, m$.\pn
In $RHS$ the general $l$ term contributes, for $i\sspeq 0,\,1,\,2$, with $\widehat{R}(l)\sspdef -R(m,\,y)\,S(l,\,y)$ to $x^{m\sspp l \sspp i}$. Therefore one  considers an $m$ rows and three columns array $AR$ with entries $AR(l,\,i)\sspeq \widehat{R}(l)\,x^{m+l+i}$. The anti-diagonals of $AR$ have identical powers of $x$.\pn
In the last row, $l \sspeq m-1$, the last two entries, and in the row $l\sspeq m-2$ the last entry are not relevant because the powers are $x^{2m}$ and $x^{2m+1}$, already treated.
\pn
A special case is $AR(0,\,0)\sspeq {\widehat R}(0)\,x^m$ because here the $LHS$ has entry $-R(m,\,y)\,x^m$, and to this power also the later discussed array $AS$ contributes with three terms $S(m,\, y), \, -y\,S(m\sspm 1,\,y)$ and $S(m\sspm 2,\,y)$ from the second terms of $AS(0,\,0)$, the first terms of  $AS(m-1,\,1)$ and  $AS(m-2,\, 2)$, respectively. Then $RHS\sspm LHS$ vanishes for $x^m$ because $-R(m,\,y) \sspp (S(m,\, y) \sspm y\,S(m\sspm 1,\,y) \sspp S(m\sspm 2,\,y)) \sspm (-R(m,\,y)) \sspeq 0$, by cancellation of $R$ and the recurrence relation of the $S$ polynomials. As announced, the $R$ term is irrelevant, only the $S$ recurrence enters.
\pn
Another case where the length of the anti-diagonal in $AR$ is not $3$ is $AR(1,\,0)\sspeq {\widehat R}(1)$ and $AR(0,\,1)\sspeq {\widehat R}(0)$. This produces for the power $x^{m+1}$  only $-R(m,\,y)(\,S(1\,,y) -y\,S(0,\,y))\sspeq 0$, because $S(-1,\,y)\sspeq 0$. There will be no contribution from the array $AS$ for this power.\pn
All other anti-diagonals of $AR$, \ie those with powers $x^{m\sspp j}$, for $j\sspin \{\range{2}{m-1}\}$, have length $3$, and their contributions $R(m,\,y)\,(S(j,\, y) \sspm y\,S(j \sspm 1,\, y)\sspp S(j \sspm 2,\, y))$ vanish because of the recurrence for the $S$-polynomials, independently of $R$.
\psn
For the pure $S$ terms of $RHS$ one considers the companion array $AS$ with two term entries $AS(l,\,i)\sspeq S(l,\,y)\,x^{l+i}\sspp S(m\sspp l,\,y)\,x^{m\sspp l\sspp i}$.
\pn 
The anti-diagonals with length 3, \ie those for powers $\{x^{j},\,x^{m\sspp j})$, for $j\sspin \{\range{2}{m-1}\}$ give vanishing contributions to $RHS$ because for both $S$-polynomial terms their recurrence relation appears.
\pn
The first term of the entry $AS(0,\,0)\sspeq S(0,\,y)\,x^0\sspp S(m,\,y)\,x^m$ is identical with $1\cdot x^0$ of $LHS$, and the second term has been needed above (among others) in the proof of the vanishing of the contribution to $x^m$.
\pn
The second anti-diagonal $A(1,\,0)$ and $A(0,\,1)$ contributes to $x^1$ with $S(1,\,y) \sspm y\,S(0,\,y)\sspeq 0$ ($S(-1,\,y)\sspeq 0$), and to $x^{m+1}$ with $S(m\sspp 1,\,y) \sspm y\,S(m,\,y)$ which is needed, together with the first term of the last entry $AS(m-1,\ 2)\sspeq S(m\sspm1,y)\,x^{m+1}$, to prove the vanishing for the contribution of $AS$ to  $x^{m+1}$. For the corresponding vanishing of the $AR$ contribution see above. There is no such $LHS$ contribution.
\psn
The second term of $AS(m-1,\,2)\sspeq S(2\,m\sspm 1,\,y)\,x^{2\,m\sspp 1}$ has been used in the treatment of the highest power above. 
\pn
Finally, the second terms of the two anti-diagonal entries $AS(m-1,\, 1)$ and $AS(m-2,\, 2)$ contribute $-y\,S(2\,m\sspm 1,\,y) \sspp  S(2\,m\sspm 2,\,y)$. They have been treated above together with the $AR$ contribution to the second highest power above.
\pn   
All $RHS$ entries of $AR$ and  $AS$ have been considered and shown to contribute only to the three powers $x^1,\, x^m$ and $x^{2\,m}$ giving the $LHS$, which ends the proof. \hskip 8.11cm $\square$
\pb
To end this work the results of the alternative bisection formulas \Eq{52} and \Eq{53} are given for the $H01$ and $H$ sequences. There derivation is done with the help of eqs. (25) and (16).\psn
\Beqarray
H01(r,s;\,2\,m\sspp 1) &\sspeq& (\sqrt{-s}\,)^{m+1}\,H01(r,s;\,m)\,R(m+1,\,r/\sqrt{-s}\,)\sspp (-s)^m\,, \nonumber \\
&\sspeq& (-s)^m\,\left(S(m-1,\,r/\sqrt{-s}\,)\,R(m+1,\,r/\sqrt{-s}\,)\sspp 1\right)\,,\\
H01(r,s;\,2\,m) &\sspeq& (\sqrt{-s}\,)^{m}\,H01(r,s;\,m)\,R(m,\,r/\sqrt{-s}\,)\,, \nonumber \\
&\sspeq& (\sqrt{-s})^{2\,m-1}\,S(m-1,\,r/\sqrt{-s}\,)\,R(m,\,r/\sqrt{-s}\,).\\
H(p,q;r,s;\,2\,m\sspp 1) &\sspeq& (\sqrt{-s}\,)^{m}\,\left\{H01(r,s;\, m)\,(\sqrt{-s}
\,q\,R(m\sspp 1,\,r/\sqrt{-s}) + p\,s\,R(m,\,r/\sqrt{-s}) ) \right. \nonumber\\
&&\hskip 1.7cm  \left. + (\sqrt{-s}\,)^m\,q   \right\}\, \nonumber\\
&\sspeq& (-s)^m\left\{S(m\sspm 1,\,r/\sqrt{-s}\,)\,(q\,R(m\sspp 1,\,r/\sqrt{-s}\,)\sspm \sqrt{-s}\,p\,R(m,\,r/\sqrt{-s}))\right. \nonumber \\ 
&& \hskip 1.3cm \left. + q \right\}\,,\\
H(p,q;r,s;\, 2\,m)  &\sspeq&  (\sqrt{-s}\,)^m\left\{R(m,\,r/\sqrt{-s}\,)\,(q\,H01(r,s;\,m) \sspp s\,p\,H01(r,s;\,m\sspm 1)) \right. \nonumber \\ 
&& \hskip 1.5cm \left.\sspm (\sqrt{-s}\,)^m\,p \right\}\,,\nonumber \\
&\sspeq&  (-s)^{m-1}\left\{R(m,\,r/\sqrt{-s}\,)\,(\sqrt{-s}\, q\,S(m\sspm 1,\,r/\sqrt{-s}\,) \sspp s\,p\,S(m\sspm 2,\,r/\sqrt{-s}\,)) \right. \nonumber\\
&& \hskip 2cm \left. \sspp s\,p \right\}\,.
\Eeqarray

\pbn
\hrule   
\pbn
OEIS\cite{OEIS} A-numbers:\seqnum{A0000045},\,\seqnum{A01445},\,\seqnum{A015448},\,\seqnum{A033887},\,\seqnum{A034807},\,\seqnum{A049310},\,\seqnum{A087960},\,\seqnum{A127672}.
\
\psn
Keywords: multi-section of sequences, Vandermonde, three-term recurrence.
\psn
MSC-Numbers:11B37, 11B39, 11B83.
   
\pbn
\hrule
\pbn
\eject

\end{document}